\documentclass[final]{amsart}
\usepackage{amsfonts,amssymb}
\usepackage[notref,notcite]{showkeys}

\usepackage[ps2pdf]{hyperref}

\newcommand{\snseqns}{\eqref{e:X-def}--\eqref{e:X-idata}}
\newcommand{\slburgers}{\eqref{e:burgers-xdef}--\eqref{e:burgers-xidata}}

\newcommand{\del}{\partial}
\newcommand{\lap}{\triangle}
\newcommand{\inv}{^{-1}}
\newcommand{\transpose}{^\text{t}}
\newcommand{\leqs}{\leqslant}
\newcommand{\geqs}{\geqslant}
\newcommand{\grad}{\nabla}
\newcommand{\gradt}{\grad\transpose}
\newcommand{\divergence}{\grad \cdot}
\newcommand{\curl}{\grad \times}

\newcommand{\womega}{{w^\omega}}
\newcommand{\urandom}{\tilde{u}}

\newcommand{\E}{\boldsymbol{\mathrm{E}}}
\newcommand{\lhp}{\boldsymbol{\mathrm{P}}}

\newcommand{\at}[1]{\big|_{#1}}

\newcommand{\D}[2][u]{\del_t #2 + \left(#1 \cdot \grad\right) #2 - \nu \lap #2}
\newcommand{\Dzero}[2][u]{\del_t #2 + \left(#1 \cdot \grad\right) #2}
\newcommand{\Dop}{\mathcal{D}}

\newcommand{\Ifn}{I}
\newcommand{\R}{\mathbb{R}}
\newcommand{\holderspace}[2]{\ensuremath{C^{\ifx0#1{#2}\else{#1,#2}\fi}}}

\newif\iftextstyle
\textstyletrue
\everydisplay\expandafter{\the\everydisplay\textstylefalse}

\newcommand{\hnorm}[3]{\iftextstyle\|#1\|\else\left\|#1\right\|\fi_{\ifx0#2{#3}\else{#2,#3}\fi}}

\newcommand{\qv}[2]{\iftextstyle\langle#1,#2\rangle\else\left\langle#1\,,\,#2\right\rangle\fi}

\numberwithin{equation}{section}
\allowdisplaybreaks

\newtheorem{theorem}{Theorem}[section]
\newtheorem{lemma}[theorem]{Lemma}
\newtheorem{prop}[theorem]{Proposition}
\newtheorem{cor}[theorem]{Corollary}

\theoremstyle{definition}
\newtheorem{defn}[theorem]{Definition}

\theoremstyle{remark}
\newtheorem{rem}[theorem]{Remark}
\newtheorem*{rem*}{Remark}

\begin{document}
\title[A stochastic representation of the Navier-Stokes equations]{A stochastic Lagrangian representation of the $3$-dimensional incompressible Navier-Stokes equations}
\author{Peter Constantin}
\thanks{The first author's work is partially supported by NSF-DMS grant 0504213.}
\address{%
Department of Mathematics\\
University of Chicago\\
Chicago IL 60637}
\email{const@math.uchicago.edu}
\author{Gautam Iyer}
\address{%
Department of Mathematics\\
University of Chicago\\
Chicago IL 60637}
\email{gautam@math.uchicago.edu}
\keywords{stochastic Lagrangian, incompressible Navier-Stokes}
\subjclass[2000]{Primary 76D05, 
60K40
}
\begin{abstract}
In this paper we derive a probabilistic representation of the deterministic $3$-dimensional Navier-Stokes equations based on stochastic Lagrangian paths. The particle trajectories obey SDEs driven by a uniform Wiener process; the inviscid Weber formula for the Euler equations of ideal fluids is used to recover the velocity field. This method admits a self-contained proof of local existence for the nonlinear stochastic system, and can be extended to formulate stochastic representations of related hydrodynamic-type equations, including viscous Burgers equations and LANS-alpha models.
\end{abstract}
\maketitle

\section{Introduction}\label{s:intro}
%

Formalisms representing solutions of partial differential equations as the expected value of functionals of stochastic processes date back to Einstein and Feynman in physics, Kolmogorov and Kac in mathematics. The theory for linear parabolic equations is well developed \cite{friedman,ks}. The theory for nonlinear partial differential equations is not as well developed and it  involves either branching processes \cite{legall}, or implicit, fixed point representations. Our representation is of this latter kind, `nonlinear in the sense of McKean'
: the drift term in the stochastic differential equation is computed implicitly
 as the expected value of an expression involving the flow it drives.
 
Connections between stochastic evolution and the deterministic Navier-Stokes equations have been established in seminal work of Chorin \cite{chorin}. In two dimensions, the nonlinear equation obeyed by the vorticity has the form a Fokker-Planck (forward Kolmogorov) equation. The Biot-Savart law relates the vorticity to the velocity in a linear fashion. These facts are used by Chorin to formulate the random vortex method to represent the vorticity of the Navier-Stokes equation using random walks and a particle limit. The random vortex method has been proved to converge by Goodman \cite{goodman} and Long \cite{long}, see also \cite{mb}. A stochastic representation of the Navier-Stokes equations for two dimensional flows using similar ideas but without 
discretization is given in \cite{bf2d}. 

A heuristic representation of the Navier-Stokes equations in three
dimensions using ideas of random walks and particles limits was done by Peskin \cite{peskin}. There are many examples of modeling approaches using stochastic 
representations in physical situations, for instance \cite{chorin1, pope}. 

The three dimensional situation is complicated by the fact that there are no obvious Fokker-Planck like equations describing the solutions. LeJan and Sznitman \cite{sznit} used a backward-in-time branching process in Fourier space to express the velocity field of a three dimensional viscous fluid as the average of a stochastic process. Their approach did not involve a limiting process, and this led to a new existence theorem. This was later generalized \cite{ossiander} to a physical space analogue. 

More recently, Busnello, Flandoli and Romito developed a representation for the $3$-dimensional Navier-Stokes equations using noisy flow paths, similar to the ones considered here. They used Girsanov transformations to recover the velocity field of the fluid from magnetization variables, and generalized their method to work for any Fokker-Planck-type equation. Their formulation however does not admit a self contained local existence proof, which is possible using our formulation.

\smallskip

In this work we describe a  stochastic representation of deterministic hydrodynamic equations that is based on the combination of the active vector formulation of the inviscid equation and the shift by a uniform Wiener process. The procedure consists thus of two steps. In the first step a Weber formula is used to express the velocity field of the inviscid equation in terms of the particle trajectories of the inviscid equation \textit{without} involving time derivatives. The second step is to replace classical Lagrangian trajectories by stochastic flows driven by the velocity field. Averaging the stochastic trajectories produces the solution. It is essential that time derivatives do not appear in the expressions to be averaged.


Perhaps the simplest illustration of this approach is the stochastic formulation of the viscous Burgers equation. The inviscid Burgers equation is
\begin{equation*}
\Dzero{u} = 0
\end{equation*}
with initial data
\begin{equation*}
u(x,0) = u_0(x).
\end{equation*}
The absence of a pressure term allows the velocity to be transported by the fluid flow. Thus if $X$ is the fluid flow map, $u \circ X$ is constant in time and $u = u_0 \circ (X_t\inv)$ gives a Weber formula to recover the velocity from the instantaneous flow map. More explicitly, the system
\begin{align*}
\dot{X} &= u\\
A_t &= X_t\inv\\
u &= u_0(A)
\intertext{with initial data}
X(a,0) &= a
\end{align*}
is equivalent to the inviscid Burgers equation before the formation of shocks. By the notation $X_t \inv$, we mean the spatial inverse of the map $X$ at time $t$.

It can be shown that if we add uniform noise to particle trajectories, then the average of the stochastic system satisfies the viscous Burgers equation. (We remark that we do not have to restrict our attention to gradient solutions $u = \grad \phi$.)
\begin{theorem}\label{t:sl-burgers}
Let $W$ be a $n$-dimensional Wiener process, $k \geqs 1$ and $u_0 \in C^{k+1, \alpha}$. Let the pair $u, X$ be a solution of the stochastic system
\begin{align}
\label{e:burgers-xdef} dX &= u\,dt + \sqrt{2 \nu}\,dW\\
A &= X\inv\\
\label{e:burgers-udef} u &= \E [u_0 \circ A]
\intertext{with initial data}
\label{e:burgers-xidata} X(a, 0) &= a.
\end{align}
For boundary conditions, we demand that either $u$ and $X - I$ are spatially periodic, or that $u$ and $X-I$ decay at infinity. Then $u$ satisfies the viscous Burgers equations
\begin{equation}
\label{e:burger-uevol} \D{u} = 0
\end{equation}
with initial data $u_0$. Here $\E$ denotes the expected value with respect to the Wiener measure and $I$ is the identity function.
\end{theorem}

\begin{rem}\label{r:burgers-invert}
The spatial invertibility of $X$ can be seen as follows: An elementary computation shows
$$ \det (\grad X_t) = \exp{ \int_0^t \divergence u_s\,ds}$$
showing $X_t$ is locally invertible for all time. Global invertibility now follows since our boundary conditions ensure $X_t$ is properly homotopic to the identity map and hence has degree $1$. We also remark that since the noise $W_t$ is spatially constant, the inverse map $A_t$ is as (spatially) regular as the flow map $X_t$. We refer the reader to \cite{sperturb, thesis, kunita} for the details.
\end{rem}

\begin{rem}\label{r:burgers-longtime}
If a solution to the system \slburgers\ exists on the time interval $[0,T]$, then our proof will show that $u$ satisfies equation \eqref{e:burger-uevol} on this time interval. Though global existence for \eqref{e:burger-uevol} is known, fixed point methods will only yield a local existence result for \slburgers.

Conversely, given a (global) solution $u$ of \eqref{e:burger-uevol} which is either spatially periodic or decays at infinity, standard theory shows that \eqref{e:burgers-xdef} has a global solution. Now uniqueness of strong solutions for linear parabolic equations and Corollary \eqref{c:spde-v} shows that \eqref{e:burgers-udef} is satisfied for all time. We refer the reader to \cite{thesis} for the details.
\end{rem}

We provide two proofs of the corresponding result for the Navier-Stokes equations (Theorem \ref{t:sl-ns}), each of which can be adapted to yield the above theorem. We present one proof of Theorem \ref{t:sl-burgers} in Section \ref{s:ito-proof}, and leave the details of the second proof to the interested reader. In the next section we describe our stochastic formulation of the Navier-Stokes equations.

\section{The stochastic Lagrangian formulation of the incompressible Navier-Stokes equations}\label{s:sl-intro}

In this section we show how the method described in the introduction can be used to produce a stochastic formulation of the Navier-Stokes equations. We begin by recalling a method of recovering the velocity field from the instantaneous flow map of an inviscid fluid. This is the Eulerian-Lagrangian formulation \cite{ele} which describes the Euler equations as an active vector system.

\begin{prop}\label{p:el-euler}
Let  $k \geqs 0$ and $u_0 \in \holderspace{k+1}{\alpha}$ be divergence free. Then $u$ satisfies the incompressible Euler equations
\begin{gather}
\Dzero{u} + \grad p = 0\\
\divergence{u} = 0
\end{gather}
with initial data $u_0$ if and only if the pair of functions $u$, $X$ satisfies the system
\begin{align}
\label{e:ele-Xdef} \dot{X} &= u\\
\label{e:ele-Adef} A &= X\inv\\
\label{e:ele-udef} u &= \lhp[ (\gradt A) \, (u_0\circ A) ]
\intertext{with initial data}
\label{e:ele-Xidata} X(a,0) &= a.
\end{align}
Here $\lhp$ is the Leray-Hodge projection \cite{cm,stein} on divergence free vector fields. We impose either periodic boundary conditions and demand that $u$ and $X - \Ifn$ are spatially periodic, or demand that $u$ and $X - \Ifn$ decay sufficiently rapidly at infinity ($\Ifn$ is the identity map).
\end{prop}

The proof of one direction of this theorem follows immediately from Lemma \ref{l:w-evol}. A complete proof can be found in \cite{ele}.
Now to obtain a solution to the viscous system we use exactly the same Weber formula \eqref{e:ele-udef} for $u$, but consider noisy trajectories instead of deterministic ones and average the noise out.

\begin{theorem}\label{t:sl-ns}
Let $\nu > 0$, $W$ be an $n$-dimensional Wiener process, $k\geqs1$ and $u_0 \in \holderspace{k+1}{\alpha}$ be a given deterministic divergence free vector field. Let the pair $u$, $X$ satisfy the stochastic system
\begin{align}
\label{e:X-def} dX &= u \,dt + \sqrt{2 \nu} \,dW\\
\label{e:A-def} A &= X\inv\\
\label{e:u-def} u &= \E \lhp\left[ (\gradt A) \, (u_0 \circ A) \right]
\intertext{with initial data}
\label{e:X-idata} X(a,0) &= a.
\end{align}
We impose boundary conditions by requiring $u$ and $X-I$ are either spatially periodic, or decay sufficiently at infinity. Then $u$ satisfies the incompressible Navier-Stokes equations
\begin{gather*}
\D{u} + \grad p = 0\\
\divergence{u} = 0
\end{gather*}
with initial data $u_0$.
\end{theorem}

\begin{rem}
Remarks \ref{r:burgers-invert} and \ref{r:burgers-longtime} are also applicable here. See also \cite{thesis}.
\end{rem}

\begin{rem}\label{r:sl-ns-force}
In the presence of a deterministic external force $f$, we only need to replace $u_0$ in equation \eqref{e:u-def} with $\varphi$ defined by
$$\varphi_t = u_0 + \int_0^t (\gradt X) f(X_s, s)\, ds.$$
We provide a proof of this in Section \ref{s:ito-proof}, however we remark that the proof of Theorem \ref{t:sl-ns} given in Section \ref{s:ode-proof} can also be adapted to yield this. Clearly, if the forcing is random and independent of the Wiener process $W$, then  our procedure provides a representation of the stochastically forced Navier-Stokes equations.
\end{rem}

\begin{rem}
The construction above can be modified to provide a stochastic representation of the LANS-alpha (or Camassa-Holm) equations. The inviscid Camassa-Holm \cite{chen, holm} equations are
\begin{gather*}
\Dzero{v} + (\gradt u) v + \grad p = 0\\
u = (1 - \alpha^2\lap)\inv v\\
\divergence v = 0
\end{gather*}
Lemma \ref{l:w-evol} gives a formula to recover $v$ from the inverse of the flow map. Thus we obtain a stochastic representation of the viscous Camassa-Holm equations by replacing \eqref{e:u-def} in \snseqns\ with
\begin{gather}
\label{e:ch-vdef} v = \E\lhp\left[ (\gradt A) \, u_0 \circ A \right]\\
\label{e:ch-udef} u = (1 - \alpha^2 \lap)\inv v.
\end{gather}
The velocity $v$ will now satisfy the viscous equation
\begin{equation*}
\Dzero{v} + (\gradt u) v - \nu \lap v + \grad p = 0.
\end{equation*}
We draw attention to the fact that the diffusive term is $\nu \lap v$ and not $\nu \lap u$. However, we do not derive the relation $u= (1-\alpha^2\lap)\inv v$; and any other translation-invariant filter $u = Tv$ would work as well.
\end{rem}

We provide two independent proofs of Theorem \ref{t:sl-ns}. We postpone these proofs to sections \ref{s:ode-proof} and \ref{s:ito-proof} respectively, and devote the remainder of this section to consequences of this theorem.

The first consequence we mention is that we have a self contained proof for the local existence of the stochastic system \snseqns.

\begin{theorem}[Local existence]\label{t:snsexist}
Let $k \geqs 1$ and $u_0 \in \holderspace{k+1}{\alpha}$ be divergence free. There exists a time $T = T(k, \alpha, L, \hnorm{u_0}{k+1}{\alpha})$, but independent of viscosity, and a pair of functions $\lambda, u \in C([0,T], \holderspace{k+1}{\alpha})$ such that $u$ and $X = \Ifn + \lambda$ satisfy the system \snseqns. Further $\exists U = U( k, \alpha, L, \hnorm{u_0}{k+1}{\alpha})$ such that for all $t \in [0,T]$ we have $\hnorm{u(t)}{k+1}{\alpha} \leqs U$.
\end{theorem}

Here $L > 0$ is a given length scale, and \holderspace{k}{\alpha} is the space of $(k, \alpha)$ H\"older continuous functions which are periodic with period $L$. The theorem is still true if we consider the domain $\R^3$, and impose decay at infinity boundary conditions instead. We do not present the proof in this paper, but refer the reader to \cite{sperturb, thesis}.

We remark again that our estimates, and existence time are independent of viscosity. The theorem and proof also work when the viscosity $\nu = 0$, thus we have a proof that gives us local existence of both the Euler and Navier-Stokes equations.

The nature of our formulation causes most identities for the Euler equations (in the Eulerian-Lagrangian form) to be valid in the above stochastic formulation after averaging. We begin by presenting identities for the vorticity.

\begin{prop}\label{p:vort-rep}
Let $\omega = \curl u$ be the vorticity, and $\omega_0 = \curl u_0$ be the initial vorticity. Then
\begin{equation}
\label{e:vort-trans} \omega = \E \left[ \left( (\grad X) \,\omega_0\right) \circ A \right].
\end{equation}
If the flow is two dimensional then the above formula reduces to
\begin{equation}
\label{e:2d-vort-trans} \omega = \E \left[\omega_0 \circ A\right].
\end{equation}
\end{prop}

\begin{rem}
In the presence of an external force $f$, we have to replace $\omega_0$ in equations \eqref{e:vort-trans} and \eqref{e:2d-vort-trans} with $\varpi$ defined by
$$ \varpi_t = \omega_0 + \int_0^t (\grad X_s)\inv g(X_s, s) \, ds$$
where $g = \curl f$. For two dimensional flows this reduces to
$$ \varpi_t = \omega_0 + \int_0^t g(X_s, s) \, ds.$$
We draw attention to the fact that these are exactly the same as the expressions in the inviscid case.
\end{rem}

We can prove this proposition in two ways: The first method is to directly differentiate \eqref{e:u-def}, and use the fact that $\lhp$ vanishes on gradients. We leave the details of this to the interested reader. The other method (presented in Section \ref{s:ito-proof}) is to use Proposition \ref{p:spde-a} and the generalized It\^o formula to show that $\omega$ satisfies the vorticity equation
\begin{equation}
\label{e:vort-evol}\D{\omega} = (\omega \cdot \grad) u.
\end{equation}

The two dimensional vorticity equation does not have the stretching term $(\omega \cdot \grad) u$.  With our formulation the two dimensional vorticity equation is immediately obtained by observing that the third component of the flow map satisfies $X_3(a,t) = a_3$.

Although it is evident from equations \eqref{e:u-def} and \eqref{e:vort-trans}, we explicitly point out that the source of growth in the velocity and vorticity fields arises from the gradient of the noisy flow map $X$. The Beale-Kato-Majda \cite{bkm} criterion guarantees if the vorticity $\omega$ stays bounded, then no blow up can occur in the Euler equations.  In the case of the Navier-Stokes equations, well known criteria for regularity exist and they can be translated in criteria for the average of the stochastic flow map. As is well known, the vorticity of a two dimensional fluid stays bounded, which immediately follows from equation \eqref{e:2d-vort-trans}.

Finally we mention the conservation of circulation. For this we need to consider a stochastic velocity $\urandom$ defined by
\begin{equation}
\label{e:urandom-def}\urandom = \lhp \left[(\gradt A) \, (u_0 \circ A) \right].
\end{equation}
We remark that $\urandom$ is defined for all realizations of the Wiener process, because the diffusion matrix in equation \eqref{e:X-def} is spatially constant. This can be seen from the techniques used by LeBris, Lions \cite{lebris} or those used in Section \ref{s:ode-proof}. Notice immediately that $u = \E \urandom$ and $\urandom_0 = u_0$. The circulation of the stochastic velocity is conserved by the stochastic flow.

\begin{prop}
If $\Gamma$ is a closed curve in space, then
$$ \oint_{X(\Gamma)} \urandom \cdot dr = \oint_\Gamma u_0 \cdot dr.$$
\end{prop}
\begin{proof}
By definition of $\lhp$, there exists a function $q$ so that
\begin{alignat*}{2}
&&	\urandom &= (\gradt A) (u_0\circ A) + \grad q\\
&\implies\quad&	\gradt X\at{A} \urandom &= u_0 \circ A + \gradt X\at{A} \grad q\\
&\implies& (\gradt X) (\urandom \circ X) &= u_0 + \grad(q \circ X).
\end{alignat*}
Hence
\begin{align*}
\oint_{X(\Gamma)} \urandom \cdot dr &= \int_0^1 (\urandom\circ X \circ \Gamma) \cdot (\grad X\at{\Gamma} \Gamma') \,dt\\
    &= \int_0^1 (\gradt X\at{\Gamma}) (\urandom \circ X \circ \Gamma) \cdot \Gamma' \,dt\\
    &= \oint_\Gamma (u_0 + \grad( q \circ X) ) \cdot dr = \oint_\Gamma u_0 \cdot dr.
\end{align*}

We remark that the above proof is exactly the same as a proof showing circulation is conserved in inviscid flows.
\end{proof}

The rest of this paper is devoted to proving the results stated in Sections \ref{s:intro} and \ref{s:sl-intro}. We provide two independent proofs of Theorem \ref{t:sl-ns}, however only adapt the second proof to yield Theorems \ref{t:sl-burgers}, \ref{t:sl-ns} and Proposition \ref{p:vort-rep}.

\section{Proof of the stochastic representation using pointwise solutions}\label{s:ode-proof}

In this section we prove Theorem \ref{t:sl-ns} by constructing pointwise (in the probability space) solutions to the SDE \eqref{e:X-def}. This idea has been used by LeBris and Lions in \cite{lebris} using a generalization of the $W^{1,1}$ theory. In our context however, the velocity $u$ is spatially regular enough for us to explicitly construct the pointwise solution without appealing to the generalized $W^{1,1}$ theory.

We remark that the proof given here will also prove Theorem \ref{t:sl-burgers}, and leave the details to the interested reader. We begin with a few preliminaries.

\begin{defn}
Given a (divergence free) velocity $u$, we define the operator $\Dop$ by
\begin{equation*}
\Dop v = \Dzero{v}
\end{equation*}
\end{defn}

\begin{lemma}
The commutator $[\Dop, \grad]$ is given by
\begin{equation*}
[\Dop, \grad] f = \Dop (\grad f) - \grad( \Dop f) = -\gradt u \grad f
\end{equation*}
\end{lemma}

\begin{proof}
By definition,
\begin{align*}
[\Dop, \grad] f &= \Dop (\grad f) - \grad( \Dop f) \\
    &= (u \cdot \grad) \grad f - \grad\left[ (u \cdot \grad) f \right]\\
    &= (u \cdot \grad) \grad f - (u \cdot \grad) \grad f - \gradt u \grad f \qedhere
\end{align*}
\end{proof}

\begin{lemma}\label{l:w-evol}
Given a divergence-free velocity $u$, let $X$ and $A$ be defined by
\begin{gather*}
\dot{X} = u(X)\\
X(a,0) = a\\
A = X\inv
\end{gather*}
We define $v$ by the evolution equation
\begin{equation*}
\Dop v = \Gamma.
\end{equation*}
with initial data $v_0$. If $w$ is defined by
\begin{equation*}
w = \lhp\left[ \left(\gradt A\right) \: v\right],
\end{equation*}
then the evolution of $w$ is given by the system
\begin{gather*}
\del_t w + (u \cdot \grad) w + (\gradt u) w + \grad p = (\gradt A) \Gamma\\
\divergence w = 0\\
w_0 = \lhp v_0
\end{gather*}
\end{lemma}

\begin{proof}
By definition of the Leray-Hodge projection, there exists a function $p$ such that
\begin{align*}
w &= \gradt A \: v - \grad p\\
&= v_i \grad A_i - \grad p\\
\implies\quad \Dop w &= (\Dop v_i) \grad A_i + v_i \Dop \grad A_i - \Dop \grad p\\
&= \Gamma_i \grad A_i - v_i (\gradt u) \grad A_i - \grad \Dop p + (\gradt u) \grad p\\
&= (\gradt A) \Gamma - (\gradt u)( v_i \grad A_i + \grad p) - \grad \Dop p\\
&= (\gradt A) \Gamma - (\gradt u)w - \grad \Dop p. \qedhere
\end{align*}
\end{proof}

\begin{cor}
If $u, X, A$ are as above, and we define $w$ by
\begin{equation*}
w = \lhp\left[ (\gradt A) \: u_0 \circ A \right].
\end{equation*}
Then $w$ evolves according to
\begin{gather*}
\Dop w + (\gradt u) w + \grad p = 0\\
\divergence w = 0\\
w(x,0) = w_0(x)
\end{gather*}
\end{cor}

\begin{proof}
The proof follows from Lemma \ref{l:w-evol} by setting $v=u_0 \circ A$ and $\Gamma = 0$.
\end{proof}

We now return to the proof of Theorem \ref{t:sl-ns}.
\begin{proof}[Proof of Theorem \ref{t:sl-ns}]

For simplicity and without loss of generality we take $\nu = \frac{1}{2}$. Let $(\Omega, \mathcal{F}, P)$ be a probability space and $W:[0,\infty) \times \Omega \to \R^3$ a three dimensional Wiener process. Define $u^\omega$ and $Y^\omega$ by
\begin{equation*}
u^\omega(x,t) = u(x + W_t(\omega), t)
\end{equation*}
and
\begin{gather*}
\dot{Y}^\omega = u^\omega(Y^\omega)\\
Y^\omega(a,0) = a
\end{gather*}
Although $u^\omega$ is not Lipschitz in time, it is certainly uniformly (in time) Lipschitz in space. Thus the regular Picard iteration will produce solutions of this equations. Finally notice that the map $X$ defined by
\begin{equation*}
X(a, t, \omega) = Y^\omega(a, t, \omega) + W_t(\omega)
\end{equation*}
solves the SDE \eqref{e:X-def}.

Let $B^\omega$ be the spatial inverse of $Y^\omega$. Notice that
\begin{align*}
X(B^\omega(x - W_t, t), t) &= Y^\omega( B^\omega(x - W_t, t), t) + W_t\\
    &= x
\end{align*}
and hence
\begin{equation}
A = \tau_{W_t} B^\omega
\end{equation}
where $\tau_x$ is the translation operator defined by
\begin{equation}
\tau_x f (y) = f(y - x)
\end{equation}

We define $\womega$ by
\begin{equation*}
\womega = \lhp\left[ (\gradt B^\omega) \; u_0 \circ B^\omega \right].
\end{equation*}
By Lemma \ref{l:w-evol},  the function $\womega$ evolves according to
\begin{gather}
\label{e:womega-evol} \del_t \womega + (u^\omega \cdot \grad) \womega + (\gradt u^\omega) w + \grad q^\omega = 0\\
\label{e:womega-dfree} \divergence \womega = 0\\
\label{e:womega-idata} \womega(x,0) = u_0(x).
\end{gather}

Now using equation \eqref{e:u-def} we have
\begin{alignat}{2}
\nonumber &&	u &= \E \lhp \left[ (\gradt A) \; u_0 \circ A \right]\\
\nonumber &&	&= \E \lhp \left[ (\gradt \tau_{W_t} B^\omega) \; u_0 \circ \tau_{W_t} B^\omega \right]\\
\nonumber &&	&= \E \lhp \left[ \tau_{W_t} \left( (\gradt B^\omega) \; u_0 \circ B^\omega \right)\right]\\
\nonumber &&	&= \E \tau_{W_t} \lhp \left[ (\gradt B^\omega) \; u_0 \circ B^\omega \right]\\
\label{e:u=ew} &&	&= \E \tau_{W_t} \womega
\end{alignat}

Our assumption $u_0 \in \holderspace{k+1}{\alpha}$ along with Theorem \ref{t:snsexist} guarantee that $\womega$ is spatially regular enough to apply the generalized It\^o formula \cite{kunita} to $\womega( x - W_t, t)$, and we have
\begin{multline*}
\womega( x - W_t, t) - u_0(x) = \int_0^t \womega( x - W_s, ds) - \int_0^t \grad \womega \at{x-W_s, s} \, dW_s + \\
    + \tfrac{1}{2} \int_0^t \lap \womega\at{x-W_s, s} \,ds + \qv{\int_0^t \del_j \womega( x - W_s, ds)}{x_j - W^j_t}.
\end{multline*}
Notice that the process $\womega$ is $C^1$ in time (since the time derivative is given by equation \eqref{e:womega-evol}), and hence bounded variation. Thus the joint quadratic variation term vanishes. Taking expected values and using \eqref{e:u=ew} we conclude
\begin{equation}
\label{e:u1} u(x,t) - u_0(x) = \E \int_0^t \womega( x - W_s, ds) + \tfrac{1}{2} \int_0^t \lap u(x,s) \,ds
\end{equation}

Using equation \eqref{e:womega-evol} and the definition of the It\^o integral we have
\begin{align}
\nonumber \E \int_0^t \womega( x - W_s, ds) &= \E \int_0^t \del_t \womega\at{x-W_s, s} ds\\
\nonumber     & = - \E \int_0^t \left[ ( u^\omega \cdot \grad ) \womega + (\gradt u^\omega) \womega + \grad q^\omega \right]_{x-W_s, s} ds\\
\nonumber     & = - \E \int_0^t \Big[ \left( u(x,s) \cdot \grad \right) \womega\at{x-W_s, s} + \left(\gradt u(x,s)\right) \womega\at{x-W_s, s} +\\
\nonumber     &\qquad\qquad + \grad q^\omega\at{x-W_s, s}\Big]\, ds\\
\nonumber     & = - \int_0^t \Big[ \left( u(x,s) \cdot \grad \right) u\at{x, s} + \left(\gradt u(x,t)\right) u\at{x, s} + \\
\nonumber     & \qquad\qquad + \grad \E q^\omega(x - W_s, s) \Big] \, ds\\
\label{e:Eintwomega}     & = - \int_0^t \left[ \left( u(x,s) \cdot \grad \right) u\at{x, s} + \grad q'\at{x,s} \right]\,ds
\end{align}
where $q'$ is defined by
\begin{equation*}
q' = \tfrac{1}{2} \grad |u|^2 + \E \tau_{W_t} q^\omega
\end{equation*}

Using equations \eqref{e:Eintwomega} in \eqref{e:u1} (along with the observation that the joint quadratic variation term is $0$) we obtain
\begin{equation*}
u(x,t) - u_0(x) = - \int_0^t \left[ (u \cdot \grad) u + \tfrac{1}{2} \lap u + \grad q \right]_{x,s} \, ds
\end{equation*}
Equations \eqref{e:womega-dfree} and \eqref{e:u=ew} show that $u$ is divergence free, concluding the proof.
\end{proof}

\section{Proof of the stochastic representation using the It\^o formula}\label{s:ito-proof}

In this section we provide a proof of Theorems \ref{t:sl-burgers}, \ref{t:sl-ns} and Proposition \ref{p:vort-rep} directly using the generalized It\^o formula. At the end of the section we provide a brief contrast with the diffusive Lagrangian formulation \cite{elns}.\smallskip

Before beginning our computations we remark that our assumptions $k\geqs 1$, $u_0 \in \holderspace{k+1}{\alpha}$ and Theorem~\ref{t:snsexist} guarantee and that the processes $A$ and $X$ are spatially regular enough to apply the generalized It\^o formula~\cite{kunita}. Further, as shown in \cite{sperturb} the displacements $X-\Ifn$ and $A-\Ifn$ are spatially $\holderspace{k+1}{\alpha}$ (thus bounded in periodic domains) and hence the It\^o integrals (and expectations) that arise here are all well defined. This said, we assume subsequently that all processes are spatially regular enough for our computations to be valid.

We begin with a lemma leading up to computing the It\^o derivative of $A$.
\begin{lemma}\label{l:mpart-a}
There exists a process $B$ of bounded variation such that almost surely
\begin{equation}
\label{e:mpart-a} A_t = B_t - \sqrt{2\nu}\int_0^t \grad A_s \,dW_s
\end{equation}
\end{lemma}
\begin{proof}
We apply the generalized It\^o formula to $A \circ X$ to obtain (almost surely)
\begin{align}
\nonumber 0 &= \int_{t'}^t A(X_s, ds) + \int_{t'}^t \grad A\at{X_s,s} \,dX_s + \tfrac{1}{2}\int_{t'}^t \del^2_{ij} A\at{X_s, s} d\qv{X_i}{X_j}_s + \\
\nonumber     &\qquad\qquad + \qv{\int_{t'}^t \del_i A(X_s, ds)}{X^i_t - X^i_{t'}}\\
\label{e:daox} &= \int_{t'}^t A(X_s, ds) + \int_{t'}^t \left[\grad A\at{X_s,s} u + \nu\lap A\at{X_s, s} \right]ds + \sqrt{2\nu}\int_{t'}^t \grad A\at{X_s,s} \,dW_s + \\
\nonumber    &\qquad\qquad + \qv{\int_{t'}^t \del_i A(X_s, ds)}{X^i_t - X^i_{t'}}.
\end{align}
Notice that the second and fourth terms on the right are of bounded variation. Since the above equality holds for all $t', t$, and $X_s$ is a homeomorphism, the lemma follows.
\end{proof}

We remark that the $\lap A$ term in \eqref{e:daox} has a positive sign, which is the anti-diffusive sign forward in time. The quadratic variation term in \eqref{e:daox} brings in twice the negative Laplacian, thus correcting this sign, and we end up with  a dissipative SPDE for $A$ as we should.

\begin{prop}\label{p:spde-a}
The process $A$ satisfies the stochastic partial differential equation
\begin{equation}
\label{e:spde-a} dA_t + (u \cdot \grad) A \, dt - \nu \lap A \,dt + \sqrt{2\nu} \grad A \, dW_t = 0
\end{equation}
\end{prop}
\begin{proof}
We apply the generalized It\^o formula to $A\circ X$ as in Lemma \ref{l:mpart-a}. Notice that differentiating equation \eqref{e:mpart-a}, we find the martingale part of $\del_i A$. Since the joint quadratic variation term in equation \eqref{e:daox} depends only on the martingale part of $\del_i A$, we can now compute it explicitly.
\begin{align*}
\qv{\int_{t'}^t \del_i A(X_s, ds)}{X^i_t - X^i_{t'}}
    &= -2\nu \qv{\int_{t'}^t \grad \left[\del_i A\right]_{X_s, s} dW_s }{W^i_t - W^i_{t'}}\\
    &= -2\nu \int_{t'}^t \lap A\at{X_s, s} ds	\qquad\text{a.s.}
\end{align*}
Substituting this in equation \eqref{e:daox}, the proposition follows.
\end{proof}

\begin{cor}\label{c:spde-v}
Let $\vartheta$ be spatially $C^2$, and differentiable in time. Then the process $v = \vartheta \circ A$ satisfies the stochastic PDE
\begin{equation}
\label{e:spde-v} dv_t + (u \cdot \grad) v_t\, dt -\nu \lap v_t \,dt + \sqrt{2\nu} \grad v_t \,dW_t = \del_t \vartheta\at{A_t} \,dt
\end{equation}
\end{cor}

\begin{rem}\label{r:traj-const}
Consider the case when $\vartheta$ is independent of time. Now when $\nu = 0$, we know that equation \eqref{e:spde-a} is equivalent to the fact that $A$ is constant along particle trajectories. Thus any function of $A$, in particular $v$, is also constant along trajectores and hence automatically satisfies equation \eqref{e:spde-v}. When $\nu \neq 0$, Proposition \ref{p:spde-a} and Corolary \ref{c:spde-v} make the same assertion for noisy trajectories.
\end{rem}
\begin{proof}
The corollary follows directly from Proposition \ref{p:spde-a} and the generalized It\^o formula:
\begin{align*}
dv_t &= \del_t \vartheta\at{A_t}\,dt + \grad \vartheta\at{A_t} \,dA_t + \tfrac{1}{2} \del^2_{ij} \vartheta\at{A_t} d\qv{A_i}{A_j}_t\\
    &= \left[ \del_t \vartheta\at{A_t}-\grad \vartheta\at{A_t} (\grad A_t) u_t + \nu \grad \vartheta\at{A_t} \lap A_t + \nu \del^2_{ij} \vartheta\at{A_t} \del_k A^i_t \del_k A^j_t \right] dt -\\
    &\qquad\qquad - \sqrt{2\nu} \grad \vartheta\at{A_t} \grad A_t \,dW_t\\
    &= \left[ \del_t \vartheta\at{A_t}-(u_t \cdot \grad) v_t + \nu \lap v_t \right] dt- \sqrt{2\nu} \grad v_t \,dW_t\qedhere
\end{align*}
\end{proof}

The results from sections \ref{s:intro} and \ref{s:sl-intro} now follow as direct consequences of Proposition \ref{p:spde-a} and Corollary \ref{c:spde-v}.
\begin{proof}[Proof of Theorem \ref{t:sl-burgers}]
The proof of this theorem follows by setting $\vartheta = u_0$, integrating \eqref{e:spde-v} and taking expected values.
\end{proof}

\begin{proof}[Proof of Theorem \ref{t:sl-ns}]
We provide the proof in the presence of an external force $f$, as stated in Remark \ref{r:sl-ns-force}. Let $v = \varphi \circ A$, $w = (\gradt A) v$ and we compute the It\^o derivative of $w$ using It\^o's formula:
\begin{align*}
dw_i &= (\del_i A) \cdot dv + d( \del_i A) \cdot v + d\qv{\del_i A_j}{v_j}\\
    &= \del_i A \cdot \left[ (-u\cdot \grad) v + \nu \lap v + (\gradt X)\at{A} f\right] dt - \sqrt{2 \nu} \del_i A \cdot ( \grad v \, dW) +\\
    &\quad+ v \cdot \left[ -\left((\del_i u) \cdot \grad\right) A - ( u \cdot \grad) \del_i A + \nu \lap \del_i A\right] dt - \sqrt{2\nu} v \cdot ( \grad \del_i A \,dW)\\
    &\quad+ 2 \nu \del^2_{ki} A_j \del_k v_j \,dt.
\end{align*}
Making use of the identities
\begin{gather*}
(u \cdot \grad) w_i = \del_i A \cdot [ (u \cdot \grad) v] + [(u \cdot \grad) \del_i A] \cdot v\\
\lap w_i = \del_i A \cdot \lap v + \lap \del_i A \cdot v + 2 \del_{ki} A_j \del_k v_j\\
\del_i u_k w_k = v \cdot \left[ \left( \del_i u \cdot \grad\right) A \right]\\
\del_i A \cdot \left[ (\gradt X)\at{A} f\right] = f
\end{gather*}
and
$$\del_k w_i = \del_i A_j \del_k v_j + v_j \del_{ki} A_j$$
we conclude
\begin{equation}
\label{e:dw} dw = \left[ -(u \cdot \grad) w + \nu \lap w - (\gradt u) w + f \right] dt - \sqrt{2\nu} \grad w \,dW.
\end{equation}
Now from equation \eqref{e:u-def} we have (almost surely)
\begin{alignat*}{2}
&&	u &= \E w + \grad q\\
&\implies& u - u_0 &= \E \int_0^t \left[ -(u \cdot \grad) w + \nu \lap w - (\gradt u) w + f\right] + \grad q\\
&&	&= \int_0^t \left[ -(u \cdot \grad) (u - \grad q) + \nu \lap (u - \grad q) - (\gradt u) (u - \grad q) + f\right] + \grad q\\
&&	&= \int_0^t \left[ -(u \cdot \grad)u + \nu \lap u + f \right] + \grad p
\end{alignat*}
where
$$p = q - \int_0^t \left[ (u \cdot \grad) q - \nu \lap q + \tfrac{1}{2} |u|^2 \right].$$

Differentiating immediately yields the theorem.
\end{proof}

\begin{proof}[Proof of Proposition \ref{p:vort-rep}]
Notice first that $\grad X$ is differentiable in time. We set $\vartheta = (\grad X) \omega_0$, $\tilde{\omega} = \vartheta \circ A$ and apply Corollary \ref{c:spde-v} to obtain
\begin{align*}
d\tilde{\omega} + (u \cdot \grad) \tilde{\omega} \, dt - \nu \lap \tilde{\omega}\, dt + \sqrt{2\nu} \grad \tilde{\omega}\, dW &= \grad \del_t X\at{A_t} \omega_0(A) \,dt\\
    &= (\grad u) (\grad X)\at A \omega_0(A) \, dt\\
    &= (\grad u) \tilde{\omega} \,dt
\end{align*}
Integrating and taking expected values shows that $\omega = \E \tilde{\omega}$ satisfies the vorticity equation \eqref{e:vort-evol} with initial data $\omega_0$. The proposition follows now follows from the uniqueness of strong solutions.
\end{proof}
\subsection{A comparison with the diffusive Lagrangian formulation.}
The computations above illustrate the connection between the stochastic Lagrangian formulation \snseqns\ presented here, and the deterministic diffusive Lagrangian formulation. We briefly discuss this below. In \cite{elns} the Navier-Stokes equations were shown to be equivalent to the system
\begin{gather*}
\D{\bar{A}} = 0\\
u = \lhp[ (\gradt \bar{A}) v]\\
\D{\bar{v}_\beta} = 2 \nu C^i_{j,\beta} \del_j \bar{v}_i\\
C^p_{j,i} = (\grad \bar{A})\inv_{ki} \del_k \del_j \bar{A}_p
\end{gather*}
with initial data
\begin{gather*}
\bar{A}(x,0) = 0\\
\bar{v}(x,0) = u_0(x).
\end{gather*}

We see first that $\bar{A} = \E A$. The commutator coefficients $C^\alpha_{ij}$ in the evolution of $\bar{v}$ compensate for the first order terms in $\lap ( (\gradt \bar{A}) \bar{v} )$. In the stochastic formulation these arises naturally from the generalized It\^o formula as the joint quadratic variation term $\qv{ \del_i A_j}{ v_j}$. More explicitly, the equations
\begin{alignat*}{2}
&&	2\nu (\gradt \bar{A}) C^j_{k, \cdot} \del_k \bar{v}_j &= 2\nu \del^2_{kj} \bar{A}_i \del_k \bar{v}_j\\
&\text{and}\qquad& d\qv{\del_i A_j}{ v_j} &= 2\nu \del^2_{kj} A_i \del_k v_j dt
\end{alignat*}
illustrate the connection between the two representations.

\section*{Acknowledgment}
Stimulating discussions with L. Ryzhik and E. S. Titi are gratefully acknowledged. P.C.'s research was partially supported by NSF-DMS grant 0504213.

\end{document}